\documentstyle[12pt]{article}

\font\twlgot =eufm10 scaled \magstep1
\font\egtgot =eufm8
\font\sevgot =eufm7
\font\twlmsb =msbm10 scaled \magstep1
\font\egtmsb =msbm8
\font\sevmsb =msbm7

\newfam\gotfam

\textfont\gotfam\twlgot
\scriptfont\gotfam\egtgot
\scriptscriptfont\gotfam\sevgot

\newfam\msbfam
\textfont\msbfam\twlmsb
\scriptfont\msbfam\egtmsb
\scriptscriptfont\msbfam\sevmsb
\def\Bbb{\protect\pBbb}
\def\pBbb{\relax\ifmmode\expandafter\Bb\else\typeout{You cann't use
Bbb in text mode}\fi}
\def\Bb #1{{\fam\msbfam\relax#1}}

\def\thebibliography#1{\section*{References}\list
   {[\arabic{enumi}]}{\settowidth\labelwidth{#1}\leftmargin\labelwidth
     \advance\leftmargin\labelsep
     \usecounter{enumi}}
     \def\newblock{\hskip .11em plus .33em minus .07em}
     \sloppy\clubpenalty4000\widowpenalty4000
     \sfcode`\.=1000\relax}

\def\op#1{\mathop{\fam0 #1}\limits}

\newcommand{\beq}{\begin{equation}}
\newcommand{\eeq}{\end{equation}}
\newcommand{\ben}{\begin{eqnarray}}
\newcommand{\een}{\end{eqnarray}}
\newcommand{\be}{\begin{eqnarray*}}
\newcommand{\ee}{\end{eqnarray*}}
\newcommand{\bea}{\begin{eqalph}}
\newcommand{\eea}{\end{eqalph}}

\newcommand{\cH}{{\cal H}}
\newcommand{\cF}{{\cal F}}

\newcommand{\bL}{{\bf L}}

\newcommand{\al}{\alpha}

\newcommand{\bt}{\beta}
\newcommand{\dl}{\delta}
\newcommand{\la}{\lambda}

\newcommand{\f}{\phi}

\newcommand{\Om}{\Omega}
\newcommand{\m}{\mu}

\newcommand{\g}{\gamma}

\newcommand{\vt}{\vartheta}
\newcommand{\vf}{\varphi}

\newcommand{\w}{\wedge}

\newcommand{\ol}{\overline}
\newcommand{\dr}{\partial}
\newcommand{\ar}{\op\longrightarrow}

\newcommand{\ve}{\varepsilon}

\newcounter{eqalph}
\newcounter{equationa}
\newcounter{theorem}
\newcounter{remark}
\newcounter{proposition}
\newcounter{lemma}
\newcounter{corollary}
\newcounter{definition}

\newenvironment{eqalph}{\stepcounter{equation}
\setcounter{equationa}{\value{equation}}
\setcounter{equation}{0}

\begin{eqnarray}}{\end{eqnarray}\setcounter{equation}{\value{equationa}}}

\def\theremark{\arabic{remark}}

\def\thedefinition{\arabic{definition}}

\newenvironment{proof}{\noindent
{\it Proof.}}{\medskip}

\newenvironment{theo}{\refstepcounter{definition}
\bigskip\noindent{\it Theorem \thedefinition.}}{\medskip}

\newcommand{\mar}[1]{}

\begin{document}
\hbox{}

{\parindent=0pt

{\large\bf Action-angle coordinates for time-dependent completely
\\ integrable Hamiltonian systems}
\bigskip

{\sc Giovanni Giachetta\footnote{E-mail address:
giovanni.giachetta@unicam.it}, Luigi
Mangiarotti$\dagger$\footnote{E-mail
address: luigi.mangiarotti@unicam.it} and Gennadi
Sardanashvily$\ddagger$\footnote{E-mail address: sard@grav.phys.msu.su;
URL: http://webcenter.ru/$\sim$sardan/ }}
\bigskip

\begin{small}
$\dagger$ Department of Mathematics and Physics, University of Camerino, 62032
Camerino (MC), Italy \\
$\ddagger$ Department of Theoretical Physics, Physics Faculty, Moscow State
University, 117234 Moscow, Russia
\bigskip



{\bf Abstract.}
A time-dependent completely integrable
Hamiltonian system is proved to admit the
action-angle coordinates around any instantly compact regular 
invariant manifold.
Written relative to these coordinates, its Hamiltonian and first integrals
are functions only of action coordinates.

\end{small}
}

\section{Introduction}

A time-dependent Hamiltonian system of $m$ degrees of freedom is called
a completely integrable system (henceforth CIS)
if it admits $m$ independent first integrals
in involution. In order to provide it with action-angle coordinates, we
use the fact that a time-dependent CIS of $m$ degrees of
freedom can be extended to an
autonomous Hamiltonian system of $m+1$ degrees of
freedom where the time is regarded as a dynamic variable \cite{bouq,dew,lich}.
We show that it is an autonomous CIS. 
However, the
classical theorem \cite{arn,laz} on action-angle coordinates 
can not be applied to this CIS 
since its invariant manifolds
are never compact because of the time axis.
Generalizing this theorem, we first prove that
there is a system of action-angle coordinates on an open neighbourhood
$U$ of a regular
invariant manifold $M$ of an autonomous CIS
if Hamiltonian vector fields of first integrals on $U$ are complete
and the foliation of $U$ by invariant manifolds is trivial.
If $M$ is compact, these conditions always hold
\cite{laz}. Afterwards, we show that, if a regular connected invariant manifold
of a time-dependent CIS
is compact at each instant, it
is diffeomorphic to the product
of the time
axis $\Bbb R$ and an $m$-dimensional torus $T^m$, and it admits
an open neighbourhood equipped with
the time-dependent action-angle coordinates ,
$i=1,\ldots,m$, where $t$ is the Cartesian coordinate on $\Bbb R$ and
$\f^i$ are cyclic coordinates on $T^m$.
Written with respect to these coordinates, a Hamiltonian
and first integrals of a time-dependent CIS are functions only
of action coordinates $I_i$.

For instance, there are action-angle coordinates $(\ol I_i;\ol\f^i)$ such
that a Hamiltonian 
of a time-dependent CIS vanishes. They are paticular initial data
coordinates, constant along trajectories of a Hamiltonian system. 
Furthermore, given an arbitrary smooth function $\cH$ on $\Bbb R^m$,
there exist action-angle coordinates $(I_i;\f^i)$, obtained by the relevant
time-dependent
canonical transformations of $(\ol I_i;\ol\f^i)$, such that 
a Hamiltonian of a time-dependent CIS with respect to these coordinates
equals $\cH(I_i)$.
Thus, time-dependent action-angle coordinates 
provide a solution of the problem of representing a
Hamiltonian of a time-dependent CIS in terms of first integrals
\cite{kaus,lew}. However, this 
representation need not hold with respect to any coordinate system
because a
Hamiltonian fails to be a scalar under time-dependent
canonical transformations.

\section{Time-dependent completely integrable Hamiltonian systems}

Recall that the configuration space of a time-dependent mechanical system
is a fibre bundle $Q\to \Bbb R$
over the time axis $\Bbb R$ equipped with the bundle
coordinates $(t,q^k)$, $k=1,\ldots,m$.
The corresponding momentum phase space is the vertical
cotangent bundle
$V^*Q$ of $Q\to\Bbb R$ endowed with holonomic
coordinates $(t,q^k,p_k=\dot q_k)$ [8-10].
The cotangent bundle
$T^*Q$, coordinated by $(q^\la,p_\la)=(t,q^k,p_0,p_k)$,
is the homogeneous momentum phase space of
time-dependent mechanics. It is provided
with the canonical Liouville form $\Xi=p_\la dq^\la$, the
canonical symplectic form
$\Om=dp_\la\w dq^\la$, and the corresponding Poisson bracket
\mar{z7}\beq
\{f,f'\}_T =\dr^\la f\dr_\la f' -\dr_\la
f\dr^\la f', \qquad f,f'\in C^\infty(T^*Q). \label{z7}
\eeq
There is the one-dimensional trivial affine bundle
\mar{z11}\beq
\zeta:T^*Q\to V^*Q. \label{z11}
\eeq
Given its global section $h$, one can equip $T^*Q$ with
the global fibre coordinate $r=p_0-h$.
The fibre bundle (\ref{z11}) provides
the vertical cotangent bundle $V^*Q$ with the canonical Poisson
structure $\{,\}_V$ such that
\mar{m72',2}\ben
&& \zeta^*\{f,f'\}_V=\{\zeta^*f,\zeta^*f'\}_T, \qquad \forall
f,f'\in C^\infty(V^*Q), \label{m72'} \\
&& \{f,f'\}_V = \dr^kf\dr_kf'-\dr_kf\dr^kf'. \label{m72}
\een

A Hamiltonian of time-dependent mechanics is defined
as a global section
\be
h:V^*Q\to T^*Q, \qquad p_0\circ
h=-\cH(t,q^j,p_j),
\ee
of the affine bundle $\zeta$ (\ref{z11}) \cite{book98,sard98}. It
yields the pull-back Hamiltonian form
\mar{b4210}\beq
H=h^*\Xi= p_k dq^k -\cH dt  \label{b4210}
\eeq
on $V^*Q$. Then there exists a unique
vector field $\g_H$ on $V^*Q$ such that
\mar{z3}\ben
&& \g_H\rfloor dt=1, \qquad \g_H\rfloor dH=0, \nonumber \\
&& \g_H=\dr_t + \dr^k\cH\dr_k- \dr_k\cH\dr^k. \label{z3}
\een
Its trajectories obey the Hamilton equation
\mar{z20}\beq
\dot q^k=\dr^k\cH, \qquad \dot p_k=-\dr_k\cH. \label{z20}
\eeq

A first integral of
the Hamilton equation (\ref{z20}) is a smooth real function $F$
on $V^*Q$ whose Lie derivative
\be
\bL_{\g_H} F=\g_H\rfloor dF=\dr_tF +\{\cH,F\}_V
\ee
along the vector field $\g_H$ (\ref{z3}) vanishes, i.e., $F$
is constant on trajectories of $\g_H$.
A time-dependent Hamiltonian system $(V^*Q,H)$ is said to be
completely integrable
if the Hamilton equation (\ref{z20}) admits $m$ first integrals
$F_k$ which are
in involution with respect to the Poisson bracket $\{,\}_V$ (\ref{m72}),
and whose differentials $dF_k$ are linearly independent almost everywhere
(i.e., the set of points where this condition fails is nowhere dense).
One can associate to this CIS an autonomous CIS on $T^*Q$ as follows.

Let us consider the pull-back
$\zeta^*H$
of the
Hamiltonian form
$H$ (\ref{b4210}) onto the cotangent bundle $T^*Q$. It is readily observed that
\mar{mm16}\beq
\cH^*=\dr_t\rfloor(\Xi-\zeta^* h^*\Xi)=p_0+\cH \label{mm16}
\eeq
is a function on $T^*Q$.
Let us regard $\cH^*$
as a Hamiltonian of an autonomous Hamiltonian system on the symplectic
manifold $(T^*Q,\Om)$ \cite{jmp00}. Its Hamiltonian vector field
\mar{z5}\beq
\g_T=\dr_t -\dr_t\cH\dr^0+ \dr^k\cH\dr_k- \dr_k\cH\dr^k \label{z5}
\eeq
is projected onto the vector field $\g_H$ (\ref{z3}) on $V^*Q$ so that
\be
\zeta^*(\bL_{\g_H}f)=\{\cH^*,\zeta^*f\}_T, \qquad
\forall f\in C^\infty(V^*Q).
\ee
An immediate consequence of this relation is the following.

\begin{theo} \label{z6} \mar{z6}
(i) Given a time-dependent CIS $(\cH; F_k)$ on $V^*Q$, the
Hamiltonian system
$(\cH^*,\zeta^*F_k)$ on $T^*Q$ is a CIS.
(ii) Let $N$ be a connected regular invariant manifold of $(\cH; F_k)$.
Then $h(N)\subset
T^*Q$ is a connected regular invariant manifold of
the autonomous CIS $(\cH^*,\zeta^*F_k)$.
\end{theo}

Hereafter, we assume that the
vector field $\g_H$ (\ref{z3}) is complete. In this case, the Hamilton equation
(\ref{z20}) admits a unique global solution through each point of the
momentum phase space $V^*Q$, and trajectories of $\g_H$
define a trivial fibre bundle $V^*Q\to V^*_tQ$
over any fibre
$V^*_tQ$ of $V^*Q\to \Bbb R$. Without loss of generality, we choose
the fibre $i_0:V^*_0Q\to V^*Q$
at $t=0$.  Since $N$ is an
invariant manifold, the fibration
\mar{ww}\beq
\xi:V^*Q\to V^*_0Q \label{ww}
\eeq
also yields the fibration of $N$ onto
$N_0=N\cap V^*_0Q$
such that $N\cong \Bbb R\times N_0$
is a trivial bundle.

\section{Time-dependent action-angle coordinates}

Let us introduce the
action-angle coordinates around an
invariant manifold $N$ of a time-dependent CIS on $V^*Q$ by use of 
the action-angle
coordinates around the invariant manifold $h(N)$ of the
autonomous CIS on $T^*Q$ in
Theorem \ref{z6}. Since $N$ and, consequently, $h(N)$ are non-compact, we
first prove the following.

\begin{theo} \label{z8} \mar{z8}
Let $M$ be a connected invariant manifold of an autonomous CIS
  $\{F_\la\}$, $\la=1,\ldots,n$, on a symplectic manifold
$(Z,\Om_Z)$. Let $U$ be an open neighbourhood of $M$ such that: (i) the
differentials $dF_\la$ are independent everywhere on $U$, (ii) the
Hamiltonian vector fields $\vt_\la$ of the first integrals 
$F_\la$ on $U$ is complete, and (iii) the submersion $\times F_\la:
U\to \Bbb R^n$ 
is a trivial bundle of invariant manifolds
over a domain $V'\subset \Bbb R^n$.
Then $U$ is isomorphic
to the symplectic annulus
\mar{z10}\beq
W'=V'\times(\Bbb R^{n-m}\times T^m), \label{z10}
\eeq
provided with the action-angle coordinates
\mar{z11'}\beq
(I_1,\ldots,I_n; x^1,\ldots, x^{n-m}; \f^1,\ldots,\f^m) \label{z11'}
\eeq
such that the symplectic form on $W'$ reads
\be
\Om_Z=dI_a\w dx^a +dI_i\w d\f^i,
\ee
and the first integrals $F_\la$ depend only on
the action coordinates $I_\al$.
\end{theo}

\begin{proof}
In accordance with the well-known theorem \cite{arn},
the invariant
manifold $M$ is diffeomorphic to the product $\Bbb R^{n-m}\times T^m$,
which is the group space of the quotient $G=\Bbb R^n/\Bbb Z^m$
of the group $\Bbb R^n$ generated by Hamiltonian vector fields
$\vt_\la$ of first integrals $F_\la$ on $M$.
Namely, $M$
is provided with the group space coordinates $(y^\la)=(s^a,\vf^i)$
where $\vf^i$ are linear functions
of parameters $s^\la$ along integral curves of the
Hamiltonian vector fields $\vt_\la$ on $U$. Let
$(J_\la)$ be coordinates on
$V'$ which
are values of first integrals $F_\la$.
Let us choose a trivialization
of the fibre bundle $U\to V$ seen as a principal bundle with the
structure group $G$. We fix its global section $\chi$.
Since parameters $s^\la$ are given up to a shift, let us
provide each fibre $M_J$, $J\in V$, with the group space
coordinates $(y^\la)$
centred at the point
$\chi(J)$. Then
$(J_\la;y^\la)$ are bundle coordinates
on the annulus $W'$ (\ref{z10}).
Since $M_J$ are Lagrangian manifolds,
the symplectic form $\Om_Z$ on $W'$
is given relative to the
bundle coordinates $(J_\la;y^\la)$ by the expression
\mar{ac1}\beq
\Om_Z=\Om^{\al\bt}dJ_\al\w dJ_\bt + \Om^\al_\bt dJ_\al\w dy^\bt. \label{ac1}
\eeq
  By the very definition of coordinates $(y^\la)$, the
Hamiltonian vector fields $\vt_\la$ of first integrals take the
coordinate form $\vt_\la=\vt_\la^\al(J_\m)\dr_\al$. Moreover, since
the cyclic group $S^1$ can not act transitively on $\Bbb R$, we have
\mar{ww25}\beq
\vt_a=\dr_a +\vt_a^i(J_\la)\dr_i, \qquad \vt_i=\vt_i^k(J_\la)\dr_k. 
\label{ww25}
\eeq
The Hamiltonian vector fields $\vt_\la$ obey the relations
\mar{ww22}\beq
\vt_\la\rfloor\Om_Z=-dJ_\la,\qquad
\Om^\al_\bt \vt^\bt_\la=\dl^\al_\la. \label{ww22}
\eeq
It follows that $\Om^\al_\bt$ is a non-degenerate matrix and
$\vt^\al_\la=(\Om^{-1})^\al_\la$, i.e., the functions $\Om^\al_\bt$
depend only on coordinates $J_\la$.
A substitution of (\ref{ww25}) into (\ref{ww22}) results in the equalities
\mar{ww30,1}\ben
&& \Om^a_b=\dl^a_b, \qquad \vt_a^\la\Om^i_\la=0, \label{ww30}\\
&& \vt^k_i\Om^j_k=\dl^j_i, \qquad \vt^k_i\Om^a_k=0. \label{ww31}
\een
The first of the equalities (\ref{ww31}) shows that the matrix $\Om^j_k$
is non-degenerate, and so is the matrix $\vt^k_i$. Then the second
one gives $\Om^a_k=0$.
By virtue of
the well-known K\"unneth formula for the de Rham cohomology of a 
product of manifolds,
the closed form $\Om_Z$ (\ref{ac1}) on $W'$ (\ref{z10})
is exact, i.e., $\Om_Z=d\Xi$ where $\Xi$  reads
\be
\Xi=\Xi^\al(J_\la,y^\la)dJ_\al + \Xi_i(J_\la) d\vf^i +
\dr_\al\Phi(J_\la,y^\la)dy^\al,
\ee
where $\Phi$ is a function on $W'$.
Taken up to an exact form, $\Xi$ is brought into the form
\mar{ac2}\beq
\Xi=\Xi'^\al(J_\la,y^\la)dJ_\al + \Xi_i(J_\la) d\vf^i. \label{ac2}
\eeq
Owing to the fact that
components of $d\Xi=\Om_Z$ are independent of $y^\la$ and obey the equalities
(\ref{ww30}) -- (\ref{ww31}), we obtain the following.

(i) $\Om^a_i=-\dr_i\Xi'^a +\dr^a\Xi_i=0$. It follows that $\dr_i\Xi'^a$ is
independent of $\vf^i$, i.e., $\Xi'^a$ is affine in $\vf^i$ and, consequently,
is independent of $\vf^i$ since $\vf^i$ are cyclic coordinate. Hence,
$\dr^a\Xi_i=0$, i.e., $\Xi_i$ is a function only of coordinates $J_j$.

(ii) $\Om^k_i=-\dr_i\Xi'^k +\dr^k\Xi_i$. Similarly to item (i), one 
shows that $\Xi'^k$
is independent of $\vf^i$ and $\Om^k_i=\dr^k\Xi_i$, i.e.,
$\dr^k\Xi_i$ is a non-degenerate matrix.

(iii) $\Om^a_b=-\dr_b\Xi'^a=\dl^a_b$. Hence, $\Xi'^a=-s^a+D^a(J_\la)$.

(iv) $\Om^i_b=-\dr_b\Xi'^i$, i.e., $\Xi'^i$ is affine in $s^a$.

In view of items (i) -- (iv), the Liouville form $\Xi$ (\ref{ac2}) reads
\be
\Xi=x^adJ_a + [D^i(J_\la) + B^i_a(J_\la)s^a]dJ_i + \Xi_i(J_j) d\vf^i,
\ee
where we put
\mar{ee1}\beq
x^a=-\Xi'^a=s^a-D^a(J_\la). \label{ee1}
\eeq
Since the matrix $\dr^k\Xi_i$ is non-degenerate,
one can introduce new coordinates $I_i=\Xi_i(J_j)$, $I_a=J_a$. Then we have
\be
\Xi=-x^adI_a + [D'^i(I_\la) + B'^i_a(I_\la)s^a]dI_i + I_i d\vf^i.
\ee
Finally, put
\mar{ee2}\beq
\f^i=\vf^i-[D'^i(I_\la) + B'^i_a(I_\la)s^a] \label{ee2}
\eeq
in order to obtain the desired action-angle coordinates 
\be
I_a=J_a, \qquad I_j(J_j), \qquad x^a=s^a +S^a(J_\la), \qquad \f^i=\vf^i+
S^i(J_\la,s^b).
\ee
These are bundle coordinates on $U\to V'$ where the
coordinate shifts (\ref{ee1}) -- (\ref{ee2}) correspond to a choice of another
trivialization of $U\to V'$.
\end{proof}

Of course, the action-angle coordinates (\ref{z11'}) by no
means are unique.
For instance, let $\cF_a$, $a=1,\ldots, n-m$ be an arbitrary smooth function
on $\Bbb R^m$.
Let us consider the canonical coordinate transformation
\mar{ww26}\beq
I'_a=I_a-\cF_a(I_j), \qquad I'_k=I_k, \qquad x'^a=x^a, \qquad \f'^i= \f^i +
x^a\dr^i\cF_a(I_j). \label{ww26}
\eeq
Then $(I'_a,I'_k; x'^a,\f'^k)$ are action-angle
coordinates on the symplectic annulus which differs from $W'$ (\ref{z10})
in another trivialization.

Now, we apply Theorem \ref{z8} to the CISs in Theorem \ref{z6}.

\begin{theo} \label{z13} \mar{z13}
Let $N$ be a connected regular invariant manifold of a time-dependent
CIS on $V^*Q$, and let the image $N_0$ of its
projection $\xi$ (\ref{ww}) be compact.
Then the invariant manifold $h(N)$ of the autonomous CIS
on $T^*Q$ has an open
neighbourhood $U$ obeying the condition of Theorem \ref{z8}.
\end{theo}

\begin{proof} (i) We first show that
functions $i_0^*F_k$ make up a CIS on the symplectic
leaf $(V^*_0Q,\Om_0)$ and
$N_0$ is its invariant manifold without critical points
(i.e., where first integrals fail to be independent).
Clearly,
the functions $i_0^*F_k$ are in involution, and $N_0$ is their
connected invariant manifold. Let us show that
the set of critical points of $\{i_0^*F_k\}$ is nowhere
dense in $V^*_0Q$ and $N_0$ has none of these points.
Let $V^*_0Q$ be equipped with some coordinates $(\ol q^k,\ol p_k)$.
Then the trivial bundle $\xi$ (\ref{ww}) is provided with the bundle
coordinates $(t,\ol q^k,\ol p_k)$ which play a role of
the initial date coordinates on the momentum phase space $V^*Q$.
Written with respect to these coordinates, the first integrals
$F_k$ become time-independent. It follows that
\mar{ww12}\beq
dF_k(y)=di_0^*F_k(\xi(y)) \label{ww12}
\eeq
for any point
$y\in V^*Q$. In particular, if $y_0\in V^*_0Q$ is a critical point
of $\{i_0^*F_k\}$, then the trajectory $\xi^{-1}(y_0)$ is a critical set
for the first integrals $\{F_k\}$. The desired statement at once
follows from this result.

(ii) Since $N_0$ obeys the condition in item (i),
there is an open neighbourhood of $N_0$
in $V^*_0Q$ isomorphic to $V\times N_0$ where $V\subset \Bbb R^m$ is
a domain, and $\{v\}\times N_0$, $v\in V$, are also invariant
manifolds in $V^*_0Q$ \cite{laz}. Then
\mar{z70}\beq
W=\xi^{-1}(V\times N_0) \cong V\times N \label{z70}
\eeq
is an open neighbourhood in
$V^*Q$ of the invariant manifold $N$ foliated by
invariant manifolds $\xi^{-1}(\{v\}\times N_0)$,
$v\in V$, of the time-dependent CIS on $V^*Q$. By virtue of
the equality (\ref{ww12}), the first integrals $\{F_k\}$ have no
critical points in $W$.
For any real number $r\in(-\ve,\ve)$, let us consider a section
\be
h_r:V^*Q\to T^*Q, \qquad p_0\circ
h_r=-\cH(t,q^j,p_j) +r,
\ee
of the affine bundle $\zeta$ (\ref{z11}). Then the images
$h_r(W)$ of
$W$ (\ref{z70}) make up an open neighbourhood $U$ of
$h(N)$ in $T^*Q$. Because $\zeta(U)=W$, the
pull-backs $\zeta^*F_k$ of first integrals $F_k$ are free from critical points
in $U$, and so is the function $\cH^*$ (\ref{mm16}). Since
the coordinate $r=p_0-h$
provides a trivialization of the affine bundle $\zeta$, the open neighbourhood
$U$ of $h(N)$ is diffeomorphic to the product
\be
(-\ve,\ve)\times h(W)\cong
(-\ve,\ve)\times V\times h(N)
\ee
which is a trivialization of the fibration
\be
\cH^*\times(\times \zeta^*F_k): U\to (-\ve,\ve)\times V.
\ee

(iii) It remains to prove that the Hamiltonian vector fields of $\cH^*$ and
$\zeta^*F_k$ on $U$ are complete. It is readily observed that the
Hamiltonian vector field $\g_T$ (\ref{z5}) of $\cH^*$ is tangent to
the manifolds $h_r(W)$, and is the image
$\g_T=Th_r\circ \g_H\circ \zeta$
of the vector field $\g_H$ (\ref{z3}).
The latter is complete on $W$, and so is $\g_T$ on
$U$. Similarly, the Hamiltonian vector field
\be
\g_k=-\dr_tF_k\dr^0 +\dr^iF_k\dr_i -\dr_iF_k\dr^i
\ee
of the function $\zeta^*F_k$ on $T^*Q$ with respect to the Poisson bracket
$\{,\}_T$ (\ref{z7}) is tangent
to the manifolds $h_r(W)$, and is the image
$\g_k=Th_r\circ \vt_k\circ \zeta$
of the Hamiltonian vector field $\vt_k$
of the first integral $F_k$ on $W$ with respect to
the Poisson bracket $\{,\}_V$ (\ref{m72}). The vector fields $\vt_k$
on $W$ are vertical relative to the fibration
$W\to\Bbb R$, and are tangent to compact manifolds.
Therefore, they are complete, and so are the vector fields $\g_k$ on $U$.
Thus, $U$ is the desired open neighbourhood
of the invariant manifold $h(N)$.
\end{proof}

In accordance with Theorem \ref{z8}, the open neighbourhood $U$ of
the invariant manifold $h(N)$ of the autonomous CIS
in Theorem \ref{z13}
is isomorphic to the symplectic annulus
\mar{z41}\beq
W'=V'\times(\Bbb R\times T^m), \qquad V'=(-\ve,\ve)\times V, \label{z41}
\eeq
provided with the action-angle coordinates
$(I_0,\ldots,I_m;t,\f^1,\ldots,\f^m)$ such that the symplectic form on
$W'$ reads
\be
\Om=dI_0\w dt +dI_k\w d\f^k
\ee
By the construction in Theorem 2, 
$I_0=J_0=\cH^*$ and the corresponding generalized angle coordinate
is $x^0=t$, while the first integrals $J_k=\zeta^*F_k$
depend only on the
action coordinates $I_i$.

Since the action coordinates $I_i$ are independent of the coordinate
$J_0$, the symplectic annulus $W'$ (\ref{z41}) inherits the fibration
\mar{z46}\beq
W'\ar^\zeta W''=V\times(\Bbb R\times T^m). \label{z46}
\eeq
By the relation similar to (\ref{m72'}), the product $W''$
(\ref{z46}), coordinated  by
$(I_i;t,\f^i)$, is provided with the Poisson structure
\be
\{f,f'\}_W = \dr^if\dr_if'-\dr_if\dr^if', \qquad f,f'\in C^\infty(W'').
\label{ww2}
\ee
Therefore, one can regard $W''$  as the momentum phase space of the
time-dependent CIS
in question around the invariant manifold $N$.

It is readily observed that the Hamiltonian vector field $\g_T$ of the
autonomous Hamiltonian
$\cH^*=I_0$ is $\g_T=\dr_t$, and so is its projection $\g_H$ (\ref{z3})
on $W''$. Consequently, the Hamilton equation (\ref{z20})
of a time-dependent CIS
with respect to the action-angle coordinates take the form
$\dot I_i=0$, $\dot\f^i=0$.
Hence, $(I_i;t,\f^i)$ are the initial data coordinates.
One can introduce such coordinates as follows. Given the fibration 
$\xi$ (\ref{ww}),
let us provide $N_0\times V\subset V^*_0Q$ in Theorem \ref{z13} with
action-angle coordinates $(\ol I_i;\ol \f^i)$ for the
CIS $\{i_0^*F_k\}$ on the symplectic leaf
$V^*_0Q$. Then, it is readily observed that $(\ol I_i;t,\ol \f^i)$ are
time-dependent action-angle coordinates
on $W''$ (\ref{z46}) such that the Hamiltonian
$\cH(\ol I_j)$ of a time-dependent CIS relative to these coordinates vanishes,
i.e., $\cH^*=\ol I_0$. Using the canonical transformations (\ref{ww26}),
one can consider time-dependent action-angle coordinates
besides the initial date ones.
Given a smooth function $\cH$ on $\Bbb R^m$,
one can provide  $W''$ with
the action-angle coordinates
\be
I_0=\ol I_0-\cH(\ol I_j), \qquad I_i=\ol I_i, \qquad \f^i=\ol \f^i +
t\dr^i\cH(\ol I_j)
\ee
such that $\cH(I_i)$ is a Hamiltonian of a time-dependent CIS on
$W''$.

\end{document}